\documentclass[12pt,Йeqno]{article}
\usepackage{cmap}
\usepackage[koi8-r]{inputenc}
\usepackage{graphicx}
\usepackage{natbib}
\usepackage{amsfonts}
\usepackage[tbtags]{amsmath}
\topskip=\baselineskip%
\newcommand{\nty}{n \to \infty}

\newcommand{\lr}{\left(}

\newcommand{\rr}{\right)}

\newcommand{\lv}{\left|}
\newcommand{\rv}{\right|}

\newcommand{\kn}{k_n}
\newcommand{\mn}{m_n}
\newcommand{\eel}{\end{lemma}}

\newcommand{\xin}{X_{i:n}}
\newcommand{\cin}{c_{i,n}}
\newcommand{\uin}{U_{i:n}}
\newcommand{\wi}{W_i}
\newcommand{\win}{W_{i:n}}
\newcommand{\xia}{\xi_{\alpha}}

\newcommand{\an}{\alpha_n}
\newcommand{\bn}{\beta_n}
\newcommand{\xib}{\xi_{1-\beta}}

\newcommand{\mau}{M_{\alpha}}
\newcommand{\eps}{\varepsilon}
\newcommand{\na}{N_{\alpha}}
\newcommand{\nb}{N_{1-\beta}}
\newcommand{\alp}{\alpha}
\newcommand{\aln}{\alpha_n}
\newcommand{\ben}{\beta_n}
\newcommand{\be}{\beta}
\newcommand{\si}{\sigma}
\newcommand{\tn}{L_{n}}
\newcommand{\bl}{\bigl(}
\newcommand{\br}{\bigr)}
\newcommand{\dn}{\delta_n}

\numberwithin{equation}{section}

\newtheorem{theorem}{Theorem}[section]

\newtheorem{lemma}{Lemma}[section]

\newtheorem{corollary}{Corollary}[section]

\newtheorem{remark}{Remark}[section]

\title{{\Huge \sf Cram\'{e}r type
moderate deviations for trimmed {\it L}-statistics}
\author{\large Nadezhda  Gribkova}
\date{
{\small \it
\centerline{St.Petersburg~State~University,~Mathematics~and~Mechanics~Faculty,}
\centerline{199034,~ Universitetskaya nab.~7/9, St.~Petersburg, Russia}
}}}
\begin{document}
\maketitle 
\noindent{\bf Abstract.}
{\small   
We establish Cram\'{e}r type  moderate deviation ({\it MD}) results  for heavy trimmed $L$-statistics; we obtain our results under a~very mild smoothness condition on the inversion $F^{-1}$ ($F$ is the underlying distribution  of i.i.d. observations) near two points, where trimming occurs, we assume also some smoothness of  weights of the $L$-statistic.  Our results complement previous work on Cram\'{e}r type  large deviations ({\it LD}) for trimmed $L$-statistics
by Gribkova~(2016) and  Callaert et~al.~(1982).
}


\section{Introduction and main results}
\label{imtro}
The theory of large deviations is one of the main branches in the probability theory and its applications.
There is an~extensive literature on this subject
for the various classes of statistics, especially for  the classical case of sums of independent random variables (see, e.g.,  \cite{petrov:1975, saul_stat:1991}) and for some types  of sums
of dependent variables, e.g., for $U$-statistics (see, e.g.,  \cite{bor_web:2003,lai_shao:2011},  and the references therein).

In contrast, there are only a~few papers on this topic for $L$-statistics.
In the case of  non-trimmed  $L$-statistics  with coefficients 
generated by a smooth on $(0,1)$ weight function, the Cram\'{e}r  type large and moderate deviations  were studied by \cite{vand_verav:1982, alesk:1991}. A~highly sharp result on Cram\'{e}r type large deviations for non-trimmed $L$-statistics with a~smooth  weight function was established by \cite{ben_zit:1990}.

For the case of heavy truncated $L$-statistics, i.e.,  the case when the weight function is zero outside some interval $[\alp,\be]\subset (0,1)$, a~result on Cram\'{e}r type large deviations was  first obtained  by \cite{call_vand_ver:1982};  more recently, the latter result was  extended and strengthened in  \cite{gri:2016}, where a~different approach than in \cite{call_vand_ver:1982} was  proposed and implemented.

To conclude this introduction  we want to mention a~paper by \cite{gao_zhao:2011}, where a~general delta method in the theory of Chernoff's type  large and moderate deviations suggested and illustrated by many examples including M-estimators and L-statistics. Some interesting results on Chernoff's type  large deviations
for (non-trimmed) $L$-statistics with smooth weight function were obtained also by \cite{boistard:2007}.

In this article we supplement our previous work  on Cram\'{e}r  type large deviations for 
trimmed $L$-statistics (cf.~\cite{gri:2016}) by some results on  moderate deviations. Our approach here is the same as  in \cite{gri:2016}:  we approximate the trimmed $L$-statistic by a~non-trimmed $L$-statistic with coefficients generated by a~smooth on $(0,1)$ weight function, where the approximating  (non-trimmed) $L$-statistic is based on order statistics corresponding to a~sample of auxiliary i.i.d.  Winsorized observations. We apply a~result on moderate deviations due to \cite{vand_verav:1982}  to the approximating  $L$-statistic and  estimate suitably the remainder term of our approximation.

Let  $X_1,X_2,\dots $ be a~sequence of independent identically
distributed (i.i.d.) real-valued random variables
(r.v.'s) with common distribution function  $F$, and for each
integer $n \geq 1$ let \ $X_{1:n}\le \dots \le X_{n:n}$ denote the
order statistics based on the sample $X_1,\dots ,X_n$.

Consider the trimmed $L$-statistic given by
\begin{equation}
\label{tn}
L_n=n^{-1}\sum_{i=\kn+1}^{n-\mn}\cin\xin,
\end{equation}
where $c_{i,n}\in \mathbb{R}$, \ $k_n$, $m_n$ are two sequences of integers such that $0\le k_n<n-m_n\le n$. Put
$\alpha_n=k_n/n$, \ $\beta_n=m_n/n$. It will be assumed throughout this paper that
\begin{equation}
\label{al_be}
\alpha_n \to \alpha, \quad \beta_n\to \beta,\quad0<\alpha<1-\beta<1,  
\end{equation}
as $n\to\infty$, i.e. we focus on the case of heavy trimmed $L$-statistic.

 Define the left-continuous inverse of $F$:
$F^{-1}(u)= \inf \{ x: F(x) \ge u \}$, \ $0<u\le 1$, \
$F^{-1}(0)=F^{-1}(0^+)$, and let $F_n$, $F_n^{-1}$ denote the
empirical distribution function and its inverse respectively.

We will consider  also the trimmed $L$-statistics with coefficients  generated by a~weight function:
 \begin{equation}
 \label{tn0}
\tn^0=n^{-1}\sum_{i=\kn+1}^{n-\mn}\cin^0\xin=\int_{\aln}^{1-\ben}J(u)F_n^{-1}(u)\,du,
\end{equation}
where $\cin^0=n\int_{(i-1)/n}^{i/n}J(u)\,du$, and  $J$ is a~function defined in an~open set $I$ such that $[\alpha,1-\beta]\subset I\subseteq(0,1)$.

To state our results, we will need the following set of assumptions.

\smallskip
\noindent{\bf (i)} {\it $J$ is Lipschitz in $I$, i.e. there exists a~constant $C\ge 0$ such that}
\begin{equation}
\label{LipJ}
|J(u)-J(v)|\le C|u-v|,\quad \forall \ \ u,\,v\in I.
\end{equation}

\smallskip
\noindent{\bf (ii)} {\it  There exists a~positive~$\eps$   such that for each $t\in\mathbb{R}$
\begin{equation}
\label{smoothF}
\begin{split}
F^{-1}\bl\alp+t\sqrt{{\log n}/{n}}\br-F^{-1}\bl\alp\br&=O\bl(\log n)^{-(1+\eps)}\br,\\
F^{-1}\bl 1-\be+t\sqrt{{\log n}/{n}}\br-F^{-1}\bl1-\be\br&=O\bl(\log n)^{-(1+\eps)}\br
\end{split}
\end{equation}
as $\nty$.}

\smallskip
\noindent{\bf (iii)} {\it $\max(|\alpha_n-\alpha|,\,|\beta_n-\beta|)=O\bl\sqrt{\frac{\log n}{n}}\br$ as $\nty$.}

\smallskip
\noindent{\bf (iv)} {\it For some~$\tilde{\eps}>0$ }
\begin{equation}
\label{c_sum}
\sum\limits_{i=k_n+1}^{n-m_n}|c_{i,n}-\cin^0|=O\bl \frac 1{\log^{\tilde{\eps}} n}\sqrt{\frac{n}{\log n}}\br.
\end{equation}

\medskip
Define a~sequence of  centering constants
 \begin{equation}
 \label{mun}
\mu_n=\int_{\aln}^{1-\ben}J(u)F^{-1}(u)\,du.
\end{equation}
Since $\aln\to\alp$,  $\ben\to\be$ as $\nty$, both variables $\tn^0$ and $\mu_n$ are well defined for all sufficiently large $n$.

It is well known  (see,~e.g.,~\cite{mas_shor:1990, stigler:1974, vandervaart:1998}) that when  the inverse $F^{-1}$ is continuous at two points $\alp$ and $1-\be$, the smoothness condition~\eqref{LipJ} implies  the weak convergence to the normal law: $\sqrt{n}(\tn^0-\mu_n)\Rightarrow N(0,\si^2)$, where 
\begin{equation}
\label{sigma}
\si^2=\si^2(J,F)=\int_{\alp}^{1-\be}\int_{\alp}^{1-\be} J(u)J(v)(u\wedge v-uv)\, dF^{-1}(u)\,dF^{-1}(v),
\end{equation}
where  $u\wedge v=\min(u,v)$; we will also use the notation $u\vee v$ for $\max(u,v)$.

Here and in the sequel, we use the convention that $\int_a^b=\int_{[a,b)}$ when integrating with respect to the left
continuous integrator $F^{-1}$. All along the article, we assume that~$\si>0$.

Define the distribution functions of the normalized $\tn$ and $\tn^0$ respectively
\begin{equation}
\label{dfs}
F_{\tn}(x) =\textbf{P}\{\sqrt{n}(\tn-\mu_n)/\si\le x\},\quad F_{\tn^0}(x) =\textbf{P}\{\sqrt{n}(\tn^0-\mu_n)/\si\le x\}.
\end{equation}

Let $\Phi$ denote the standard normal distribution function. Here is our first result on  Cram\'{e}r  type moderate  deviations for $\tn$.

\begin{theorem}
\label{thm1}
Suppose that $F^{-1}$ satisfies condition {\bf (ii)}  and that condition~{\bf (iii)} holds for the sequences $\aln$ and $\bn$. In addition, assume that there exists a~function $J$ satisfying condition {\bf (i)} such that~{\bf (iv)} holds for the weights $\cin$.
Then
\begin{equation}
\label{thm_1}
\begin{split}
1- F_{\tn}(x) &= [1-\Phi(x)](1+o(1)),\\
F_{\tn}(-x)&=\Phi(-x)(1+o(1)),
\end{split}
\end{equation}
as \,$\nty$, uniformly in the range $-A\le x\le c\sqrt{\log n}$, for each $c>0$ and  $A>0$.
\end{theorem}

The proof of our results is  relegated to Section~\ref{proof}.

Theorem~\ref{thm1} directly implies the following  corollary.
\begin{corollary}
\label{cor2}
Let $\cin=\cin^0=n\int_{(i-1)/n}^{i/n}J(u)\,du$ \,$(\kn+1 \le i\le n-m_n)$, where $J$ is a~function  satisfying
{\bf (i)}. Assume that conditions {\bf (ii)} and {\bf (iii)} are satisfied.

Then  relations~\eqref{thm_1} with  $\tn=\tn^0$ hold true, for each $c>0$ and  $A>0$, uniformly in the range $-A\le x\le c\sqrt{\log n}$.
\end{corollary}

Finally, we state  a~version of Theorem~\ref{thm1}, where the scale factor $\sigma/n^{1/2}$ is replaced by $\sqrt{{\text{\em Var}}(\tn)}$. A~very mild moment condition  will be required now to ensure the existence of the variance of $\tn$.
\begin{theorem}
\label{thm2}
Suppose that the conditions of Theorem~\ref{thm1} hold true.  In addition, assume  that \
$\textbf{E}|X_1|^{\gamma}<\infty$ for some $\gamma>0$. Then
\begin{equation}
\label{thm_2}
\frac{\sqrt{\text{\em Var}(\tn)}}{\sigma/\sqrt{n}}=1 + O\bl(\log n)^{-(1+2\nu} \br,
\end{equation}
where $\nu=\eps\wedge\tilde{\eps}$ ($\eps$, $\tilde{\eps}$ are as in~\eqref{c_sum} and \eqref{smoothF} respectively).

Moreover,  relations~\eqref{thm_1} remain valid for each $c>0$ and  $A>0$, uniformly in the range $-A\le x\le c\sqrt{\log n}$, if we replace $\sigma/n^{1/2}$ in definition of $F_{\tn}(x)$ {\em (cf.~\eqref{dfs})} by $\sqrt{{\text{\em Var}}(\tn)}$.
\end{theorem}
\section{Stochastic approximation for $\tn^0$}

Let $\xi_{\nu}=F^{-1}(\nu)$, $0<\nu<1$, be the $\nu$-th quantile of $F$ and  $\wi$ denote $X_i$ Winsorized outside of $(\xia,\xib]$. In other words
\begin{equation}
\label{2_2}
\wi=\left\{
\begin{array}{ll}
\xia,& X_i\le \xia, \\
X_i,& \xia < X_i \le \xib,\\
\xib,& \xib < X_i .
\end{array}
\right.
\end{equation}

Let $\win$ denote the order statistics, corresponding to $W_1,\dots,W_n$, the sample of $n$ i.i.d. auxiliary random variables.

Similarly as in~\cite{gri:2016}, we will approximate  $\tn$ by a~linear combination of the order statistics $\win$ with the coefficients generated by the following weight function
\begin{equation}
\label{2_3}
J_w(u)=\left\{
\begin{array}{ll}
J(\alp),& u\le \alp, \\
J(u),& \alp < u \le 1-\be,\\
J(1-\be),& 1-\be < u
\end{array}
\right.
\end{equation}
defined in $[0,1]$. It is obvious that if $J$ is Lipschitz in $I$, i.e. satisfies condition~\eqref{LipJ} with some positive constant~$C$, then the function $J_w$ is Lipschitz in $[0,1]$ with some constant $C_w\le C$.

Consider the auxiliary non-truncated $L$-statistic given by
 \begin{equation}
 \label{Ln}
\widetilde{L}_n=n^{-1}\sum_{i=1}^{n}\widetilde{c}_{i,n}\win=\int_0^1 J_w(u)G_n^{-1}(u)\,du,
\end{equation}
where $\widetilde{c}_{i,n}=n\int_{(i-1)/n}^{i/n}J_w(u)\,du$. Define the centering constants
 \begin{equation}
 \label{muL}
\mu_{\widetilde{L}_n}=\int_0^1 J_w(u)G^{-1}(u)\,du.
\end{equation}

Since $\wi$ has the finite moments of any order and  because $J_w$ is Lipschitz, the distribution of the normalized  $\widetilde{L}_n$ tends to the standard normal law (see, e.g.,~\cite{stigler:1974})
\begin{equation*}
 \label{Ln_to}
\sqrt{n}(\widetilde{L}_n-\mu_{\widetilde{L}_n})/\sigma(J_w,G) \Rightarrow N(0,1),
\end{equation*}
where the asymptotic variance  is given by
\begin{equation}
\label{sigma}
\sigma^2(J_w,G)=\int_0^1\int_0^1 J_w(u)J_w(v)(u\wedge v-uv)\, dG^{-1}(u)\,dG^{-1}(v).
\end{equation}

Observe that for $u\in(\alp,1-\be]$ we have $J_w(u)=J(u)$, $G^{-1}(u)=F^{-1}(u)$, and that $dG^{-1}(u)\equiv 0$ for $u\notin(\alp,1-\be]$. This yields the equality of the asymptotic variances
\begin{equation}
\label{sigma_eq}
\sigma^2(J_w,G)=\sigma^2(J,F)=\sigma^2
\end{equation}
of the truncated $L$-statistic $\tn^0$ and the non-truncated $L$-statistic $\widetilde{L}_n$ based on the Winsorized random variables.

Define the binomial random variable $N_{\nu}= \sharp \{i : X_i \le
\xi_{\nu} \}$, where $0<\nu <1$.
Put $A_n=\na/n$, \ $B_n=(n-\nb)/n$.

 The following lemma provides us a~useful representation which is crucial in our proofs. This lemma is proved  in  \cite[Lemma~2.1]{gri:2016}, therefore here we present  only its statement.
\begin{lemma} {\em(\cite{gri:2016})}
\label{lem_2.1}
\begin{equation}
\label{lem_2.1_1}
\tn^0-\mu_n=\widetilde{L}_n-\mu_{\widetilde{L}_n}+R_n,
\end{equation}
where $R_n=R_n^{(1)}+R_n^{(2)}$,
\begin{equation}
\label{lem_2.1_2}
R_n^{(1)}=\int_{\alp}^{A_n} J_w(u)[F_n^{-1}(u)-\xia]\,du-\int_{1-\be}^{1-B_n} J_w(u)[F_n^{-1}(u)-\xib]\,du
\end{equation}
and
\begin{equation}
\label{lem_2.1_3}
 \ \ \ R_n^{(2)}=\int_{\an}^{\alp} J(u)[F_n^{-1}(u)-F^{-1}(u)]\,du-\int_{1-\bn}^{1-\be} J(u)[F_n^{-1}(u)-F^{-1}(u)]\,du.
\end{equation}
\end{lemma}

\begin{remark} {\em
It should be noted that
the method  based on the $L$-statistic  approximation was first applied in \cite{gri:2016}; it can be viewed as a~development of the approach proposed in \cite{gri_helm:2006, gri_helm:2007, gri_helm:2014}, where the second order asymptotic properties (Berry--Esseen bounds and  one term Edgeworth type expansions) for (intermediate) trimmed means and their Studentized and bootstrapped versions were established.
In  the articles mentioned we constructed   $U$-statistic type approximations for (intermediate) trimmed means, where  we used sums of  auxiliary i.i.d. Winsorized observations as the linear  terms; in order to get the second (quadratic) $U$-statistic terms, we applied  some special Bahadur--Kiefer representations of von Mises  statistic type for (intermediate) sample quantiles (cf.~\cite{gri_helm:2012})}.
\end{remark}
\section{Proofs}
\label{proof}
\noindent{\bf Proof of Theorem~\ref{thm1}}.
Obviously it suffices to prove the first of relations~\eqref{thm_1}. Set
\begin{equation}
\label{2_1}
V_n=\tn-\tn^0=n^{-1}\sum_{i=\kn+1}^{n-\mn}(\cin -\cin^0)\xin.
\end{equation}
Lemma~\ref{lem_2.1} and relation~\eqref{2_1} yield
\begin{equation}
\label{proof_1}
\tn-\mu_n=\widetilde{L}_n-\mu_{\widetilde{L}_n}+R_n +V_n.
\end{equation}
An~application of the classical Slutsky argument to~\eqref{proof_1} gives that, for $\delta>0$, $1- F_{\tn}(x)$ is  bounded above and below by
\begin{equation}
\label{proof_2}
\textbf{P}\{\sqrt{n}(\widetilde{L}_n-\mu_{\widetilde{L}_n})/\si> x-2\delta\}+\textbf{P}\{\sqrt{n}|R_n|/\si> \delta\}+\textbf{P}\{\sqrt{n}|V_n|/\si> \delta\}
\end{equation}
and
\begin{equation}
\label{proof_3}
\textbf{P}\{\sqrt{n}(\widetilde{L}_n-\mu_{\widetilde{L}_n})/\si> x+2\delta \}-\textbf{P}\{\sqrt{n}|R_n|/\si> \delta\}-\textbf{P}\{\sqrt{n}|V_n|/\si> \delta\}
\end{equation}
respectively. Fix arbitrary $c>0$ and $A>0$.  Set $\delta=\dn=\bl\log(n+1)\br^{-1/2-\eps_1}$, where $0<\eps_1<\eps\wedge\tilde{\eps}$, and $\eps$, $\tilde{\eps}$ are as in conditions {\bf( ii)} and {\bf(iv)} respectively(cf~\eqref{smoothF}-\eqref{c_sum}).
From~\eqref{proof_2} and~\eqref{proof_3} it immediately follows  that to prove our theorem it suffices to show that
\begin{equation}
\label{proof_4}
\textbf{P}\{\sqrt{n}(\widetilde{L}_n-\mu_{\widetilde{L}_n})/\si> x \pm 2\dn\}=[1-\Phi(x)](1+o(1)) ,
\end{equation}
\begin{equation}
\label{proof_5}
 \quad \quad \ \ \,\textbf{P}\{\sqrt{n}|R_n|/\si> \dn\}=[1-\Phi(x)]o(1) ,
\end{equation}
\begin{equation}
\label{proof_6}
\quad \quad \ \ \, \textbf{P}\{\sqrt{n}|V_n|/\si> \dn\}=[1-\Phi(x)]o(1) ,
\end{equation}
uniformly in the range $-A\le x\leq c\sqrt{\log n}$.

Let us prove~\eqref{proof_4}. Observe that  $\widetilde{L}_n$ represents a~non-truncated $L$-statistic based on the sample $W_1,\dots,W_n$ of i.i.d. bounded random variables; and since its weight function $J_w$ is  Lipschitz in $[0,1]$, we can apply a~result by \cite{vand_verav:1982}. Set $\Delta=2\sup_{n\ge 1}\dn=2/(\log 2)^{1/2+\eps_1}$. Since $\textbf{E}|W_i|^p\leq M<\infty$ (for each $p>0$ with some $M>0$), by Theorem~1~({\em i})  of \cite{vand_verav:1982}
\begin{equation}
\label{proof_7}
\textbf{P}\{\sqrt{n}(\widetilde{L}_n-\mu_{\widetilde{L}_n})/\si> x \pm 2\dn \}=[1-\Phi(x \pm 2\dn)](1+o(1)) ,
\end{equation}
uniformly with respect to $x$ such that $-(A+\Delta) \le x \pm 2\dn \leq c_1\sqrt{\log n}$, where we may take  $c_1>c$. Hence~\eqref{proof_7} holds uniformly in the range $-A\le x\leq  c\sqrt{\log n}$ for all sufficiently large $n$. Further, we apply Lemma~A1 of \cite{vand_verav:1982} (in which  the required asymptotic property of~$\Phi$ is given in a~very convenient form). Since $\dn \sqrt{\log n}=o(1)$, due to that lemma we obtain that  $1-\Phi(x\pm \dn)=[1-\Phi(x)](1+o(1))$ uniformly in the range $-A\leq x\leq c \sqrt{\log n}$.
Summarizing, we find that~\eqref{proof_4} is valid, uniformly in the range required.

Let us  prove~\eqref{proof_5}. First, we argue similarly to  the corresponding place in  \cite[Theorem~1.1]{gri:2016}.   Let $I_1^{(j)}$ and $I_2^{(j)}$ denote the first and the second terms of $R^{(j)}_n$ (cf.~\eqref{lem_2.1_2}--\eqref{lem_2.1_3}) respectively, $j=1,2$. Then     $R_n=I_1^{(1)}-I_2^{(1)}+I_1^{(2)}-I_2^{(2)}$ and
\begin{equation}
\label{proof_9}
\textbf{P}\{\sqrt{n}|R_n|/\si> \dn\}\le \sum_{k=1}^2 \textbf{P}\{\sqrt{n}|I^{(1)}_k|/\si> \dn/4\}
+\sum_{k=1}^2 \textbf{P}\{\sqrt{n}|I^{(2)}_k|/\si> \dn/4\} .
\end{equation}
Notice that for $x\in [-A,c\sqrt{\log n}]$
\begin{equation}
\label{np8}
\frac 1{1-\Phi(x)}\leq \frac 1{1-\Phi(c\sqrt{\log n})} \sim c\sqrt{\log n}\, n^{c^2/2}.
\end{equation}
Hence, it suffices to show that for each positive  $C$ (in particular for $C=\si/4$),
\begin{equation}
\label{proof_10}
\textbf{P}\{\sqrt{n}|I^{(j)}_k|>C \dn \}=o\lr(\log n)^{-1/2} n^{-c^2/2} \rr, \text{\ \ $k,j=1,2,$}
\end{equation}
as $\nty$. We will prove~\eqref{proof_10} for $I^{(1)}_1$ and $I^{(2)}_1$ (the treatment of $I^{(1)}_2$ and $I^{(2)}_2$ is similar and therefore omitted).

Consider $I^{(1)}_1$. First, note that if $\alp < A_n$, then $\max_{u\in (\alp,A_n)} |F_n^{-1}(u)-\xia|=\xia -X_{[n\alp]+1:n}\le \xia -X_{[n\alp]:n}$, as $F_n^{-1}$ is monotonic. Here and in what follows $[x]$ represents the greatest integer function. Similarly we find that if $A_n \le \alp$, then $\max_{u\in (A_n,\alp)} |F_n^{-1}(u)-\xia|=X_{[n\alp]:n}-\xia$. Furthermore, by the Lipschitz condition for $J$, there exists a~positive $K$ such that $\max_{u\in [0,1]}J_w(u)\le \sup_{u\in I}J(u)\le K$. This yields
\begin{equation}
\label{proof_11}
|I^{(1)}_1| =\lv \int_{\alp}^{A_n} J_w(u)[ F_n^{-1}(u)-\xia]\, du \rv  \le K |A_n-\alp| |X_{[n\alp]:n}-\xia|.
\end{equation}

Let $U_1,\dots,U_n$ be a~sample of independent $(0,1)$-uniform distributed random variables, $\uin$ -- the corresponding order statistics. Set $\mau=\sharp \{i : U_i \le \alp \}$.  Since the joint
distribution of  $\xin$ and $\na$ coincides with the joint distribution of $F^{-1}(\uin)$ and $\mau$, $i=1,\dots,n$, we have
\begin{equation}
\label{proof_n1}
\begin{split}
\textbf{P}\{\sqrt{n}|I^{(1)}_1|>C \dn \}
&\leq \textbf{P}\{ n^{-1/2} |\mau-\alp n| |F^{-1}(U_{[n\alp]:n})-F^{-1}(\alp)| >C \dn \}\\
&\leq P_1+P_2,
\end{split}
\end{equation}
where
\begin{equation*}
\label{proof_n2}
\begin{split}
P_1=&\textbf{P}\{|\mau-\alp n| > c_1\sqrt{n\log n} \}, \\
P_2=&\textbf{P}\{ |F^{-1}(U_{[n\alp]:n})-F^{-1}(\alp)| >C (\log (n+1))^{-(1+\eps_1)} \},
\end{split}
\end{equation*}
where we choose $c_1>c$. Here and in the sequel, $C$ stands for a~positive constant not depending on $n$, which may change  its value  from line to line.

For $P_1$ by Bernstein's inequality we obtain
\begin{equation}
\label{proof_n3}
P_1 \leq 2 exp(-h_n),
\end{equation}
with
$
h_n=
\frac{c_1^2\log n}{2[1+O(\sqrt{\log n/n})]} \sim \frac{c_1^2\log n}{2}.
$
Hence $P_1=o\bl(\log n)^{-1/2} n^{-c^2/2} \br$. Next we estimate $P_2$ on the r.h.s. in~\eqref{proof_n1}.
To shorten notation, let $k_{\alp}=[n\alp]$, $p_{\alp}=\textbf{E}U_{k_{\alp}:n}=k_{\alp}/(n+1)$, and note that
$0<\alp -p_{\alp}\leq (\alp+1)/(n+1)=O(1/n)$.
Define $\mathbb{V}_n(p_{\alp})=\sqrt{n}(U_{k_{\alp}:n}-p_{\alp})$ and let $\mathcal{E}$ denote the event $\{|\mathbb{V}_n(p_{\alp})|\leq c_1\sqrt{p_{\alp}\log n} \}$, where as before $c_1$ is an~arbitrary number such that $c_1>c$. Put
$M_n=|F^{-1}(p_{\alp}+c_1\sqrt{{p_{\alp}\log n}/{n}})-F^{-1}(p_{\alp})|\vee
 |F^{-1}(p_{\alp}-c_1\sqrt{{p_{\alp}\log n}/{n}})-F^{-1}(p_{\alp})|$. Then we have
\begin{equation}
\label{proof_n5}
P_2 \leq \textbf{P}\{ M_n >C (\log (n+1))^{-(1+\eps_1)} \} +\textbf{P}\{\overline{\mathcal{E}}\}.
\end{equation}
By condition~{\bf (ii)}, and because of $\eps_1>\eps$, the first probability on the r.h.s. in~\eqref{proof_n5} is zero for all sufficiently large $n$. In order to estimate the second  probability on the r.h.s. in~\eqref{proof_n5}, we can apply Inequality~1 given in  \cite[page 453]{shor_weln:1986}. Then we obtain
\begin{equation}
\label{proof_n6}
\textbf{P}\{\overline{\mathcal{E}}\} \leq \exp\Bigl[ -c_1^2\frac{\log n}{2}\Bigr]+\exp\Bigl[ -c_1^2\frac{\log n}{2}\widetilde\psi(t_n)\Bigr],
\end{equation}
where $\widetilde\psi$ is the function defined in \cite[page 453, formula~(2)]{shor_weln:1986},
$t_n=c_1\sqrt{\frac{\log n}{n}}$. Since $t_n\to 0$ as $\nty$, hence $t_n>-1$ for all sufficiently large $n$, and by Proposition~1 in \cite[page 455, relation~(12)]{shor_weln:1986}, we find that $\widetilde\psi(t_n)\geq \frac 1{1+2t_n/3}$. This and relation~\eqref{proof_n6} together imply that
\begin{equation}
\label{proof_n7}
\textbf{P}\{\overline{\mathcal{E}}\}\leq 2\exp\Bigl[ -c_2^2\frac{\log n}{2}\Bigr] =2n^{-c_2^2/2},
\end{equation}
for each $c_2$ such that $c<c_2<c_1$ and for all sufficiently large $n$. Summarizing, we get that
$P_2=o\bl(\log n)^{-1/2} n^{-c^2/2} \br$, and the desired bound~\eqref{proof_10} for $I^{(1)}_1$ follows.

Next we prove~\eqref{proof_10} for $I^{(2)}_1$. Define a~sequence of intervals $\Gamma_n=[\alp\wedge\an,\alp\vee\an+1/n)$,  then we obtain
\begin{equation}
\label{proof_11_}
|I^{(2)}_1| =\lv \int_{\an}^{\alp}J(u) [F_n^{-1}(u)-F^{-1}(u)]\, du \rv \le K |\an-\alp| \emph{D}_n,
\end{equation}
where $\emph{D}_n= \max_{i:\,i/n\in \Gamma_n}|\xin-F^{-1}(i/n)|\vee |\xin-F^{-1}((i-1)/n)|$. By condition~{\bf (iii)}, the estimate~\eqref{proof_11_} implies that
\begin{equation}
\label{proof_n12}
\textbf{P}\{\sqrt{n}|I^{(2)}_1|>C \dn \}
\leq \textbf{P}\{ \emph{D}_{n,u} >C (\log (n+1))^{-(1+\eps_1)} \},
\end{equation}
where $\emph{D}_{n,u}=\max_{i:\,i/n\in \Gamma_n}|F^{-1}(\uin)-F^{-1}(i/n)|\vee |F^{-1}(\uin)-F^{-1}((i-1)/n)|$.
Define $\mathbb{V}_n(p_i)=\sqrt{n}(\uin-p_i)$, where
$p_i=\textbf{E}\uin=i/(n+1)$,  and let $\mathcal{E}_i$ denote the event $\{|\mathbb{V}_n(p_i)|\leq c_1\sqrt{p_i\log n} \}$, where now we take $c_1$ such that $c_1^2>c^2+1/2$. Then by condition~{\bf (ii)}, we find that
\begin{equation}
\label{proof_n13}
\textbf{P}\{\sqrt{n}|I^{(2)}_1|>C \dn \}
\leq \textbf{P}\{\bigcup_{i:\,i/n\in \Gamma_n} \overline{\mathcal{E}}_i \} \leq \sum_{i:\,i/n\in \Gamma_n}\textbf{P}(\overline{\mathcal{E}}_i).
\end{equation}
Similarly as before, using Inequality~1 from \cite{shor_weln:1986}, for each $i:\,i/n\in \Gamma_n$, we obtain that
\begin{equation*}
\label{proof_n14}
\textbf{P}\{\overline{\mathcal{E}}_i\}\leq 2\exp\Bigl[ -c_2^2\frac{\log n}{2}\Bigr] =2n^{-c_2^2/2},
\end{equation*}
with some $c_2$ such that $c_1^2>c_2^2>c^2+1/2$, and since by condition~{\bf (iii)} $\sharp \{i:\,i/n\in \Gamma_n\}=O(\sqrt{n\log n})$, it follows from~\eqref{proof_n13} that $\textbf{P}\{\sqrt{n}|I^{(2)}_1|>C \dn \}=o\bl(\log n)^{-1/2} n^{-c^2/2} \br$. This completes the proof of~\eqref{proof_10}, which implies that~\eqref{proof_5} holds true uniformly in the range $-A\le x\leq c\sqrt{\log n}$.

Let us finally prove that~\eqref{proof_6} is valid uniformly in the range $-A\le x\leq c\sqrt{\log n}$. \ By condition~{\bf (iv)}, there exists $b>0$ such that
\begin{equation*}
\label{proof_21}
\sqrt{n}|V_n|\leq b (\log n)^{-(1/2+\tilde{\eps})}\bl|X_{(\kn+1):n}|\vee |X_{(n-\mn):n}|\br,
\end{equation*}
for all sufficiently large $n$. Thus,
\begin{equation*}
\label{proof_22}
\begin{split}
\textbf{P}\lr \sqrt{n}|V_n|/\sigma >\dn \rr &\leq \textbf{P}\lr |X_{(\kn+1):n}|\vee |X_{(n-\mn):n}|>C(\log (n+1))^{\tilde{\eps}-\eps_1} \rr\\
&\leq \textbf{P}_{3}+\textbf{P}_{4},
\end{split}
\end{equation*}
where $\textbf{P}_{3}=\textbf{P}\bigl(|X_{(\kn+1):n}|>C(\log (n+1))^{\tilde{\eps}-\eps_1}\bigr) $,  \
$\textbf{P}_{4}=\textbf{P}\bigl(|X_{(n-\mn):n}|>C(\log (n+1))^{\tilde{\eps}-\eps_1}\bigr) $, and $\tilde{\eps}-\eps_1>0$ by the choice of $\eps_1$. Let us estimate $\textbf{P}_{3}$ (the treatment for $\textbf{P}_{4}$ is same and therefore omitted). We have
\begin{equation}
\label{proof_23}
\begin{split}
\textbf{P}_{3}&=\textbf{P}\lr \lv F^{-1}(U_{(\kn+1):n})\rv >C(\log (n+1))^{\tilde{\eps}-\eps_1}\rr \\
&\leq \textbf{P}\lr \lv F^{-1}(U_{(\kn+1):n})- F^{-1}(\alp)\rv +\lv F^{-1}(\alp)\rv  >C(\log (n+1))^{\tilde{\eps}-\eps_1}\rr,\\
&= \textbf{P}\lr \lv F^{-1}(U_{(\kn+1):n})- F^{-1}(p_{\an})\rv >C(\log (n+1))^{\tilde{\eps}-\eps_1} (1+o(1))\rr,
\end{split}
\end{equation}
where $p_{\an}=\textbf{E}U_{(\kn+1):n}$. Arguing similarly as when estimating $P_2$ (cf.~\eqref{proof_n5}-\eqref{proof_n7}), we find that the r.h.s. of~\eqref{proof_23} is $o\bl(\log n)^{-1/2} n^{-c^2/2} \br$. This completes the proof of~\eqref{proof_6} and the theorem. $\square$

\medskip
\noindent{\bf Proof of Theorem~\ref{thm2}}. Let us first prove relation~\eqref{thm_2}. By Lemma~\ref{lem_2.1} and relation~\eqref{proof_1}, we have
\begin{equation*}
\label{proof t21}
\text{\em Var}(\tn)=\text{\em Var}(\widetilde{L}_n)+\text{\em Var}(R_n+V_n)+2 \text{\em cov}(\widetilde{L}_n,R_n+V_n).
\end{equation*}
Since $W_i$ are bounded, the conditions  in \cite[Theorem~2\,(ii), page~431]{vand_verav:1982} are satisfied, and hence
\begin{equation*}
\label{proof t22}
\sigma^{-1} n^{1/2}\sqrt{\text{\em Var}(\widetilde{L}_n)}=1+O(n^{-1/2})
\end{equation*}

Further, we have
\begin{equation*}
\label{proof t21}
\begin{split}
n|\text{\em cov}(\widetilde{L}_n,R_n+V_n)|&\leq n[\text{\em Var}(\widetilde{L}_n)\text{\em Var}(R_n+V_n)]^{1/2}\\
&= \sigma[n\text{\em Var}(R_n+V_n)]^{1/2}(1+O(n^{-1/2})).
\end{split}
\end{equation*}
The latter three relations imply that in order to prove~\eqref{thm_2}, it suffices to show that
\begin{equation}
\label{proof t23}
n \text{\em Var}(R_n+V_n)=O\bl(\log n)^{-(1+2\nu} \br,
\end{equation}
where $\nu= \eps\wedge\tilde{\eps}$, and  $\eps$, $\tilde{\eps}$ are the constants from conditions~{\bf (ii)} and~{\bf (iv)} respectively.
We have
\begin{equation}
\label{proof t24}
n \text{\em Var}(R_n+V_n)\leq n \textbf{E}(R_n+V_n)^2\leq  5n\Bigl[\sum_{k,j=1}^{2} \textbf{E}\bigl(I_k^{(j)}\bigr)^2 + \textbf{E}V_n^2 \,\Bigr],
\end{equation}
where $I_k^{(j)}$ are as in \eqref{proof_9}-\eqref{proof_10}. We will show that
\begin{equation}
\label{proof t25}
n \textbf{E}\bigl(I_1^{(j)}\bigr)^2 =O\bl (\log n)^{-2(1+\eps)}\br, \ \  n \textbf{E}\bigl(I_2^{(j)}\bigr)^2 =O\bl(\log n)^{-(1+2\eps)} \br,\ \ j=1,2,
\end{equation}
and that
\begin{equation}
\label{proof t26}
n \textbf{E}V_n^2 =O\bl(\log n)^{-(1+2\tilde{\eps})} \br.
\end{equation}
Relations \eqref{proof t24}-\eqref{proof t26} imply the desired bound~\eqref{proof t23}.

Let us prove the first relation in~\eqref{proof t25}. We will  consider in detail only the case $k=1$ (the treatment in the case $k=2$ is same and therefore omitted). Let as before $k_{\alp}=[\alpha n]$ and $k_n=\alpha_n n$. By \eqref{proof_11} and the Schwarz inequality, we have
\begin{equation*}
\label{proof t27}
\begin{split}
\textbf{E}\bigl(I_1^{(1)}\bigr)^2 \leq & K^2 [\textbf{E}(A_n-\alp)^4\textbf{E}(X_{k_{\alp}:n}-\xia)^4]^{1/2}\\
= & K^2 n^{-2}[\textbf{E}(\na-\alp n)^4\textbf{E}(X_{k_{\alp}:n}-\xia)^4]^{1/2}.
\end{split}
\end{equation*}
By well-known formula for 4-th moments of a~binomial random variable, we have $\textbf{E}(\na-\alp n)^4=3\alp^2(1-\alp^2)n^2(1+o(1))$. Thus, there exists a~positive constant $C$ independent of $n$ such that
\begin{equation}
\label{proof t28}
n\textbf{E}\bigl(I_1^{(1)}\bigr)^2 \leq C [\textbf{E}(X_{k_{\alp}:n}-\xia)^4]^{1/2}
\end{equation}
for all sufficiently large $n$.  Fix arbitrary $c>0$, $A>0$. Let $p_{\alp}$, $\mathbb{V}_n(p_{\alp})$ and the event
$\mathcal{E}$ be as when estimating $P_2$ in the proof of Theorem~\ref{thm1}, and $c_1$ is an~arbitrary constant such that $c_1>c$. Then we can write
\begin{equation*}
\label{proof t29}
\textbf{E}(X_{k_{\alp}:n}-\xia)^4= \textbf{E}[(F^{-1}(U_{k_{\alp}:n})-F^{-1}(\alp))^4 \textbf{1}_{\mathcal{E}}]
+ \textbf{E}[(F^{-1}(U_{k_{\alp}:n})-F^{-1}(\alp))^4 \textbf{1}_{\overline{\mathcal{E}}}].
\end{equation*}
By a~well known property of the order statistics (see, e.g.,~\cite[Thoerem~1]{gri:1995}, and due to our moment assumption, $\textbf{E}|F^{-1}(U_{k_{\alp}:n}|^k$ is bounded from above for each $k>0$. Then by condition~{\bf (ii)}, the latter quantity is of the order
$O\bl (\log n)^{-4(1+\eps)} +\textbf{P}(\overline{\mathcal{E}})\br=O\bl (\log n)^{-4(1+\eps)}\br$ (cf.~\eqref{proof_n7} ). This bound and~\eqref{proof t28} together imply that $n\textbf{E}\bigl(I_1^{(1)}\bigr)^2=O\bl (\log n)^{-2(1+\eps)}\br$.

Consider $I_1^{(2)}$. By condition {\bf (iii)}, there exists $L>0$ such that
$(\an-\alp)^2\leq L\log n/n$, for all sufficiently large $n$. Then in view of~\eqref{proof_11_} we obtain
\begin{equation}
\label{proof t210}
n\textbf{E}\bigl(I_1^{(2)}\bigr)^2 \leq n K^2 (\aln-\alp)^2 \textbf{E} \emph{D}_n^2 \leq LK^2\log n\textbf{E} \emph{D}_n^2.
\end{equation}
Hence, to get the second bound in~\eqref{proof t25}, it suffices to show that
\begin{equation}
\label{proof t211}
\textbf{E}\emph{D}_n^2=O\bl (\log n)^{-2(1+\eps)}\br,
\end{equation}
and since $\emph{D}_{n,u}\stackrel{d}=\emph{D}_{n,u}$, it suffices to prove~\eqref{proof t211} for $\emph{D}_{n,u}$.
Let $p_i$, $\mathbb{V}_n(p_i)$ and the event $\mathcal{E}_i$ be as  in the proof of Theorem~\ref{thm1}, when estimating of $I_1^{(2)}$ (cf.~\eqref{proof_n12}-\eqref{proof_n13}, where we now take $c_1$ such that $c^2_1>c^2+1/2$.
Let $\textbf{1}_{\mathcal{E}_i}$ denote  the indicator of the event~$\mathcal{E}_i$
Then we have
\begin{equation}
\label{proof t212}
\textbf{E}\emph{D}_{n,u}^2 =\textbf{E}\bl \emph{D}_{n,u}^2 \textbf{1}_{\bigcap_{i:\,i/n\in \Gamma_n}\mathcal{E}_i}\br +
\textbf{E}\bl \emph{D}_{n,u}^2 \textbf{1}_{\bigcup_{i:\,i/n\in \Gamma_n}\overline{\mathcal{E}}_i}\br.
\end{equation}
By condition~{\bf (ii)} the first term on the r.h.s. in~\eqref{proof t212} is $O\bl(\log n)^{-2(1+\eps)}\br$, and since $\textbf{E}\emph{D}_{n,u}^k$ is bounded from above for each $k>0$, there exists a~positive constant $M$, not depending on $n$, such that for all sufficiently large $n$
\begin{equation}
\label{proof t212n}
\textbf{E}\bl \emph{D}_{n,u}^2 \textbf{1}_{\bigcup_{i:\,i/n\in \Gamma_n}\overline{\mathcal{E}}_i}\br\leq
M\textbf{P}\big\{\bigcup_{i:\,i/n\in \Gamma_n}\overline{\mathcal{E}}_i\big\}\leq M\sum_{i:\,i/n\in \Gamma_n}\textbf{P}\{\overline{\mathcal{E}}_i\}.
\end{equation}
Similarly to the proof of Theorem~\ref{thm1} (cf.~\eqref{proof_n13}), we find that the magnitude on the r.h.s.
in~\eqref{proof t212n} is of the order $o\bl(\log n)^{-1/2} n^{-c^2/2} \br$. Summarizing, we obtain the
validity of the second  relation in~\eqref{proof t25}.

We now turn to the proof of~\eqref{proof t26}. We have
\begin{equation*}
\begin{split}
n\textbf{E}V_n^2 &\leq n^{-1} \Bigl( \sum_{i=k_n+1}^{n-m_n}|\cin -\cin^0|\Bigr)^2\textbf{E}\bigl(X^2_{k_n+1:n}\vee X^2_{n-m_n:n}\bigr).
\end{split}
\end{equation*}
Due to our moment assumption we have $\textbf{E}\bigl(X^2_{k_n+1:n}\vee X^2_{n-m_n:n}\bigr)=O(1)$, and
by condition {\bf (iv)} we  get $\bl\sum_{i=k_n+1}^{n-m_n}|\cin -\cin^0|\br^2=O\bl n\,(\log n)^{-(1+2\tilde{\eps})}\br$. These bounds and the latter displayed estimate  yield~\eqref{proof t26}.
Thus, relation~\eqref{thm_2} is proved.

In order to complete the proof of our theorem, it remains to argue the possibility of the replacement $\sigma/n^{1/2}$  by $\sqrt{{\text{\em Var}}(\tn)}$ in~\eqref{thm_1} without affecting the result. We prove it for the first  relation in~\eqref{thm_1}, for the second one it will then follow from the first one if we replace $\cin$ by $-\cin$.

Fix arbitrary $c>0$ and $A>0$, set $\lambda_n=\sigma^{-1}n^{1/2} \sqrt{{\text{\em Var}}(\tn)}$ and write
\begin{equation}
\label{proof t217}
\frac{\textbf{P}\bigl((\tn-\mu_n)/ \sqrt{{\text{\em Var}}(\tn)} >x\bigr)}{1-\Phi(x) }=\frac {1-F_{\tn}(\lambda_nx)}{1-\Phi(\lambda_n x)} \, \,\frac{1-\Phi(\lambda_n x)}{1-\Phi(x)}.
\end{equation}
Set $B=A\sup_{n\in \mathbb{N}}\lambda_n$. Since $\lambda_n \to 1$, the  number $B$ exists. Hence by Theorem~\ref{thm1}, the first ratio on the r.h.s. in ~\eqref{proof t217} tends to $1$ as $\nty$, uniformly in $x$ such that $-B \leq \lambda_n x \leq c\sqrt{\log n}$, hence in particular  uniformly in the range $-A \leq  x \leq c\sqrt{\log n}$. Furthermore, we see that $|\lambda_n-1|^{1/2} \sqrt{\log n} \to 0$, which is due to the fact that $|\lambda_n-1|^{1/2}=O\bl(\log n)^{-(1/2+\nu)}\br$.  Hence, by Lemma A1 from \cite{vand_verav:1982},  the second ratio on the r.h.s. in ~\eqref{proof t217} also tends to~$1$, uniformly in the range $-A \leq  x \leq c\sqrt{\log n}$. The theorem is proved.  $\square$

\bibliographystyle{apalike} 
\bibliography{gribkova_ref}

\begin{thebibliography}{}

\bibitem[Aleskeviciene, 1991]{alesk:1991}
Aleskeviciene, A. (1991).
\newblock Large and moderate deviations for {L}-statistics.
\newblock {\em Lithuanian Math.~J.}, 33:145--156.

\bibitem[Bentkus and Zitikis, 1990]{ben_zit:1990}
Bentkus, V. and Zitikis, R. (1990).
\newblock Probabilities of large deviations for {L}-statistics.
\newblock {\em Lithuanian Math.~J.}, 30:215--222.

\bibitem[Boistard, 2007]{boistard:2007}
Boistard, H. (2007).
\newblock Large deviations for {L}-statistics.
\newblock {\em Statist. \& Decis.}, 25:89--125.

\bibitem[Borovskikh and Weber, 2003]{bor_web:2003}
Borovskikh, Y. and Weber, N. (2003).
\newblock Large deviations of {U}-statistics. {I-II}.
\newblock {\em Lithuanian Math.~J.}, 43:\ 11--33,\ 241--261.

\bibitem[Callaert et~al., 1982]{call_vand_ver:1982}
Callaert, H., Vandemaele, M., and Veraverbeke, N. (1982).
\newblock A~{C}ram\'{e}r type large deviations theorem for trimmed linear
  combinations of order statistics.
\newblock {\em Comm. Statist. Th. Meth.}, 11:\ 2689--2698.

\bibitem[Gao and Zhao, 2011]{gao_zhao:2011}
Gao, F. and Zhao, X. (2011).
\newblock Delta method in large deviations and moderate deviations for
  estimators.
\newblock {\em Ann. Statist.}, 39:\ 1211--1240.

\bibitem[Gribkova, 1995]{gri:1995}
Gribkova, N. (1995).
\newblock Bounds for absolute moments of order statistics.
\newblock In {\em {\em (Skorokhod,~A.V., Borovskikh,~Yu.V.~eds.)} Exploring
  Stochastic Laws: Festschrift in Honor of the 70th Birthday of Acad. V.S.
  Korolyuk}, pages 129--134. VSP, Utrecht.
\newblock Available at \ ar{X}iv:1607.08066v2[math.PR].

\bibitem[Gribkova, 2016]{gri:2016}
Gribkova, N.\, V. (2016).
\newblock Cram\'{e}r type large deviations for trimmed {L}-statistics.
\newblock {\em Probab.~Math.~Statist. {\em(to appiar)}}.
\newblock Available at ar{X}iv:1507.02403[math.PR].

\bibitem[Gribkova and Helmers, 2006]{gri_helm:2006}
Gribkova, N.\, V. and Helmers, R. (2006).
\newblock The empirical {E}dgeworth expansion for a {S}tudentized trimmed mean.
\newblock {\em Math.~Methods~Statist.}, 15(1):\ 61--87.

\bibitem[Gribkova and Helmers, 2007]{gri_helm:2007}
Gribkova, N.\, V. and Helmers, R. (2007).
\newblock On the {E}dgeworth expansion and the {M} out of~{N} bootstrap
  accuracy for a {S}tudentized trimmed mean.
\newblock {\em Math.~Methods~Statist.}, 16(2):\ 142--176.

\bibitem[Gribkova and Helmers, 2012]{gri_helm:2012}
Gribkova, N.\, V. and Helmers, R. (2012).
\newblock On a {B}ahadur--{K}iefer representation of von {M}ises statistic type
  for intermediate sample quantiles.
\newblock {\em Probab.~Math.~Statist.}, 32(2):\ 255--279.

\bibitem[Gribkova and Helmers, 2014]{gri_helm:2014}
Gribkova, N.\, V. and Helmers, R. (2014).
\newblock Second order approximations for slightly trimmed means.
\newblock {\em Theory Probab. Appl.}, 58(3):\ 383--412.

\bibitem[Lai et~al., 2011]{lai_shao:2011}
Lai, T., Shao, Q., and Wang, Q. (2011).
\newblock Cram\'{e}r type moderate deviations for {S}tudentized {U}-statistics.
\newblock {\em ESAIM: Probability and Statistics}, 15:\ 168--179.

\bibitem[Mason and Shorack, 1990]{mas_shor:1990}
Mason, D. and Shorack, G. (1990).
\newblock Necessary and sufficient conditions for asymptotic normality of
  trimmed {L}-statistics.
\newblock {\em J.~Statist.~Plan.~Inference}, 25:\ 111--139.

\bibitem[Petrov, 1975]{petrov:1975}
Petrov, V.\, V. (1975).
\newblock {\em Sums of independent random variables}.
\newblock Springer-Verlag, New York.

\bibitem[Saulis and Statulevi\v{c}ius, 1991]{saul_stat:1991}
Saulis, L. and Statulevi\v{c}ius, V. (1991).
\newblock {\em Limit theorems for large deviations}.
\newblock Kluwer Academic Publishers, Dordrecht.

\bibitem[Shorack and Wellner, 1986]{shor_weln:1986}
Shorack, G.\, R. and Wellner, J.\, A. (1986).
\newblock {\em Empirical processes with application in statistics}.
\newblock Wiley, New York.

\bibitem[Stigler, 1974]{stigler:1974}
Stigler, S.\, M. (1974).
\newblock Linear functions of order statistics with smooth weight functions.
\newblock {\em Ann. Statist.}, 2:\ 676--693.

\bibitem[van~der Vaart, 1998]{vandervaart:1998}
van~der Vaart, A. (1998).
\newblock {\em Asymptotic statistics}.
\newblock Cambridge Univ.~Press, Cambridge.

\bibitem[Vandemaele and Veraverbeke, 1982]{vand_verav:1982}
Vandemaele, M. and Veraverbeke, N. (1982).
\newblock Cram\'{e}r type large deviations for linear combinations of order
  statistics.
\newblock {\em Ann. Probab.}, 10:\ 423--434.

\end{thebibliography}

\end{document}